\documentclass[11pt]{amsart}%
\usepackage{amssymb,amsmath,amsfonts,latexsym,amsthm,geometry,graphicx}
\usepackage{amsmath}
\usepackage{amsfonts}
\usepackage{amssymb}
\usepackage{graphicx}
\usepackage{color}%
\providecommand{\U}[1]{\protect\rule{.1in}{.1in}}
\providecommand{\U}[1]{\protect\rule{.1in}{.1in}}
\providecommand{\U}[1]{\protect\rule{.1in}{.1in}}
\providecommand{\U}[1]{\protect\rule{.1in}{.1in}}
\newtheorem{theorem}{Theorem}[section]

\theoremstyle{definition}

\newtheorem{remark}[theorem]{Remark}

\begin{document}
\title[Complex Bohnenblust--Hille inequality whose monomials have indices in an arbitrary set]{Complex Bohnenblust--Hille inequality whose monomials have indices in an arbitrary set}
\author[F. C. Alves]{Fernando C. Alves}
\address{Departamento de Matem\'{a}tica\\
\indent Universidade Federal da Para\'{\i}ba\\
\indent58.051-900 - Jo\~{a}o Pessoa, Brazil}
\email{cabralalvesf@gmail.com}
\author[D. M. Serrano]{Diana Marcela Serrano-Rodr\'{\i}guez}
\address{Departamento de Matem\'{a}ticas\\
\indent Universidad Nacional de Colombia\\
\indent111321 - Bogot\'{a}, Colombia}
\email{dmserrano0@gmail.com and diserranor@unal.edu.co}
\subjclass[2010]{}

\begin{abstract}
The Bohnenblust-Hille inequality and its variants have found applications in
several areas of Mathematics and related fields. The control of the constants
for the variant for complex $m$-homogeneous polynomials is of special interest
for applications in Harmonic Analysis and Number Theory. Up to now, the best
known estimates for its constants are dominated by $\kappa\left(
1+\varepsilon\right)  ^{m}$, where $\varepsilon>0$ is arbitrary and $\kappa>0$
depends on the choice of $\varepsilon$. For the special cases in which the number of variables in each monomial is bounded by some fixed number $M$, it has been shown that the optimal constant is dominated by a constant depending solely on $M$. In this note, based on a deep result of Bayart, we prove an inequality for any subset of the indices, showing how summability of arbitrary restrictions on monomials can be related to the combinatorial dimension associated with them.      

\end{abstract}
\maketitle

\section{Introduction}

In the investigation of complex $m$-homogeneous polynomials whose monomials
have a uniform bound $M$ on the number of variables, Carando, Defant, and
Sevilla-Peris have shown (in \cite{carando}) that the optimal constants of the
Bohnenblust-Hille inequality are dominated by a polynomial on $m$ whenever $M$
is fixed. Maia, Nogueira and Pellegrino (\cite{mais}) substantially improved
these results, proving that the optimal constants have a universal bound (i.e.
they are bounded by a constant dependent solely on $M$) under the same
hypotheses. Nevertheless, there has not been any noticeable progress for these
bounds as $M$ increases with $m$. In this note, we prove an inequality for any subset of the indices. This is achieved by employing a deep tool recently proved by Bayart \cite{bayart} which involves the concept of combinatorial dimension. We now turn to the detailed presentation of what has just been outlined.

A mapping $P:c_{0}\rightarrow \mathbb{C}$ is a continuous $m$-homogeneous polynomial if there exists a continuous $m$-linear $\hat{P}:c_{0} \times ...  \times c_{0}\rightarrow \mathbb{C}$ such that for every $x \in c_{0}$ we have $P(x)=\hat{P}(x,...,x)$. For each sequence $\alpha:\mathbb{N}\rightarrow \mathbb{N}\cup\{0\}$ satisfying $|\alpha|:=\sum_{j}\alpha_j=m$, we denote by $c_{\alpha}(P)$ the coefficient  of the monomial $x_{i_1}^{\alpha_{i_1}}\cdots x_{i_m}^{\alpha_{i_m}}$. Defining the norm of $P$ by $\Vert P\Vert:=\sup_{x\in B_{c_{0}}}|P(x)|$, the Bohnenblust--Hille inequality \cite{bh} for complex $m$-homogeneous
polynomials reads as follows: there is a constant $C_{m}\geq1$ such that
\[
\left(  \sum_{\left\vert \alpha\right\vert =m}|c_{\alpha}(P)|^{\frac{2m}{m+1}%
}\right)  ^{\frac{m+1}{2m}}\leq C_{m}\Vert P\Vert
\]
for all continuous $m$-homogeneous polynomials $P:c_{0}\rightarrow\mathbb{C}$. The
exact control of the growth of the constant $C_{m}$ plays a crucial role in a
vast number of applications. In 2011, it was proved in \cite{ann} that $C_{m}$
can be chosen with exponential growth, and in 2014 (\cite{bohr}) the result
was improved as follows: for any $\varepsilon>0$ there is a constant
$\kappa>0$ such that%
\begin{equation}
C_{m}\leq\kappa\left(  1+\varepsilon\right)  ^{m}. \label{2211}%
\end{equation}

\bigskip In order to improve the estimate (\ref{2211}) to a polynomial bound
the authors of \cite{carando} have shown that for integers $m,M$ with $M\leq
m$, we have%
\begin{equation}
\left(  \sum_{\alpha\in\Delta_{M}}|c_{\alpha}(P)|^{\frac{2m}{m+1}}\right)
^{\frac{m+1}{2m}}\leq2^{\frac{M}{2}}m^{\frac{M+1}{2}}\Vert P\Vert\label{ju}%
\end{equation}
for all continuous $m$-homogeneous polynomials $P:c_{0}\rightarrow\mathbb{C}$,
where
\[
w(\alpha)=\text{card }\{j:\alpha_{j}\neq0\}
\]
and
\[
\Delta_{M}=\{\alpha:|\alpha|=m,\quad w(\alpha)\leq M\}.
\]
In other words, the constant $C_{m}$ can be taken with a polynomial bound
provided that the sum is restricted to monomials with uniformly bounded number
of variables $M$. In \cite{mais} the result of \cite{carando} was improved,
using techniques of \cite{www}, by showing that under the same hypotheses the
constant $C_{m}$ can be replaced by a universal constant depending just on $M$.

To illustrate, we provide a pair of examples. First, for $j\in\mathbb{N}$ let $\sigma_{j}:\mathbb{N\rightarrow N}$
be injections such that for each $j,k$ with $j\not =k$ one has $card\left\{
i\in\mathbb{N}:\sigma_{j}\left(  i\right)  \not =\sigma_{k}\left(  i\right)
\right\}  =\aleph_{0}$ (e.g. define $\sigma_{j}(n):=p_{j}^n$, where $p_{j}$ is the $j$-th prime number). Notice that the results of \cite{carando, mais} are
useless for a continuous $m-$homogeneous polynomial%

\[
P(x)=%
{\textstyle\sum\limits_{i=1}^{\infty}}
c_{i}x_{\sigma_{1}\left(  i\right)  }x_{\sigma_{2}\left(  i\right)  }\cdots
x_{\sigma_{m}\left(  i\right)  },
\]
because in this case $M=m$.

Second, if $\sigma_{1}, \sigma_{2}, \sigma_{3}$ are injections from $\mathbb{N}\times\mathbb{N}$ to $\mathbb{N}$, consider the following polynomial\[
P(x)=%
{\textstyle\sum\limits_{i,j,k=1}^{\infty}}
c_{ijk}x_{\sigma_{1}\left(  i,j\right)  }x_{\sigma_{2}\left(  j,k\right)  }x_{\sigma_{3}\left(  k,i\right)  }.
\] 

Although it does not make sense at this point, our result will show that the former example behaves somewhat like $M=1$ whereas the latter exhibts fractional quality as though $M=3/2$. This will follow from the fact that, using notation we introduce below, in these examples we have $dim(\Lambda)=1$ and  $dim(\Lambda)=3/2$, respectively.

We use notations and notions from combinatorial dimension, as presented in
\cite{bayart}. For $\Lambda\subset\mathbb{N}^{m}$ and $n\geq0$, define
\[
\psi_{\Lambda}\left(  n\right)  :=\max\left\{  card\left(  \left(  A_{1}%
\times\cdots\times A_{m}\right)  \cap\Lambda\right)  ;\text{ }A_{i}%
\subset\mathbb{N}\text{, }card\left(  A_{i}\right)  \leq n\right\}  .
\]
\newline The combinatorial dimension of $\Lambda$, denoted by $\dim\left(
\Lambda\right)  $, is defined as
\begin{align*}
\dim\left(  \Lambda\right)   &  :=\limsup_{n\rightarrow+\infty}\frac{\log
\psi_{\Lambda}\left(  n\right)  }{\log n}\\
&  =\inf\left\{  s>0;\text{ }\exists C>0,\text{ }\psi_{\Lambda}\left(
n\right)  \leq Cn^{s}\text{ for all }n\in\mathbb{N}\right\}  .
\end{align*}

Finally, to each monomial $x_{i_{1}}^{\alpha_{i_{1}}} \cdots x_{i_{M}}^{\alpha_{M}}$ of a given $m$-homogeneous polynomial (i.e. $i_{1}<...<i_{M}$, $1\leq M \leq m$, and $\sum_{k=1}^{M}\alpha_{i_{k}}=m$), we associate  $(i_{1},\stackrel{\alpha_{i_{1}}}{...},i_{1},...,i_{M},\stackrel{\alpha_{i_{M}}}{...},i_{M}) \in \mathbb{N}^{m}$ and denote by $\Lambda$ the set of all such indices. To keep the sequence notation for indices introduced above, we also define $\Gamma_{\Lambda}$ to be the set of all sequences $\alpha$ in $\mathbb{N}\cup\{0\}$ naturally associated to the indices by defining $\alpha(i_k):=\alpha_{i_k}$ and zero otherwise.

\begin{theorem}
\label{arrocho} For all positive integers $m$ and all sets of indices $\Lambda
\subset\mathbb{N}^{m},$ there is a constant $C_{\Lambda}$ such that
\[
\left(  \sum_{\alpha\in\Gamma_{\Lambda}}|c_{\alpha}(P)|^{\frac{2m}{m+1}}\right)
^{\frac{m+1}{2m}}\leq e^{\dim\Lambda}\left(  C_{\Lambda}%
mm!\right)  ^{\frac{\dim\Lambda}{m}}\left(  \frac{2}{\sqrt{\pi}%
}\right)  ^{\frac{\left(  m-1\right)  \dim\Lambda}{m}}\Vert P\Vert
\]
for all continuous $m$-homogeneous polynomials $P:c_{0}\rightarrow\mathbb{C}$.
\end{theorem}

\section{The proof of Theorem \ref{arrocho}}

By \ \cite[Theorem 2.1]{bayart} we know that there is a constant
$C_{\Lambda}>0$ such that
\begin{align*}
&  \left(  \sum_{\left(  i_{1},\dots,i_{m}\right)  \in\Lambda%
}\left\vert T\left(  e_{i_{1}},\dots,e_{i_{m}}\right)  \right\vert
^{\frac{2\dim\Lambda}{1+\dim\Lambda}}\right)  ^{\frac
	{1+\dim\Lambda}{2\dim\Lambda}}\\
&  \leq C_{\Lambda}\sum_{k=1}^{m}\left(  \sum_{i_{k}=1}^{\infty
}\left(  \sum_{i_{1},...,i_{k-1},i_{k+1},...,i_{m}=1}^{\infty}\left\vert
T\left(  e_{i_{1}},\dots,e_{i_{m}}\right)  \right\vert ^{2}\right)  ^{\frac
	{1}{2}}\right)  %
\end{align*}
for all continuous $m$--linear forms $T:c_{0}\times\dots\times c_{0}%
\rightarrow\mathbb{C}$. It is not difficult to prove that using the
Khinchine's inequality, for Steinhaus variables, we have
\[
\left(  \sum_{i_{k}=1}^{\infty}\left(  \sum_{i_{1},...,i_{k-1},i_{k+1}%
	,...,i_{m}=1}^{\infty}\left\vert T\left(  e_{i_{1}},\dots,e_{i_{m}}\right)
\right\vert ^{2}\right)  ^{\frac{1}{2}}\right) \leq\left(  \frac{2}%
{\sqrt{\pi}}\right)  ^{m-1}\left\Vert T\right\Vert
\]
for all $k\in\left\{  1,...,m\right\}  $ and for all continuous $m$--linear
forms $T:c_{0}\times\dots\times c_{0}\rightarrow\mathbb{C}$. In this way,
there is a constant $C_{\Lambda}>0$ such that
\[
\left(  \sum_{\left(  i_{1},\dots,i_{m}\right)  \in\Lambda}\left\vert
T\left(  e_{i_{1}},\dots,e_{i_{m}}\right)  \right\vert ^{\frac{2\dim
		\Lambda}{1+\dim\Lambda}}\right)  ^{\frac{1+\dim
		\Lambda}{2\dim\Lambda}}\leq m\left(  \frac{2}{\sqrt{\pi}%
}\right)  ^{m-1}C_{\Lambda}\left\Vert T\right\Vert
\]
for all continuous $m$--linear forms $T:c_{0}\times\dots\times c_{0}%
\rightarrow\mathbb{C}$.

Let $\hat{P}$ be the symmetric $m$-linear form associated to $P$. Note that%
\[
\sum_{\alpha\in\Gamma_{\Lambda}}|c_{\alpha}(P)|^{\frac{2\dim\Lambda}%
	{1+\dim\Lambda}}\leq\sum_{\left(  i_{1},\dots,i_{m}\right)
	\in\Lambda}\left(  m!\right)  ^{\frac{2\dim\Lambda}%
	{1+\dim\Lambda}}\left\vert \hat{P}(e_{i_{1}},\dots,e_{i_{m}%
})\right\vert ^{\frac{2\dim\Lambda}{1+\dim\Lambda}}.
\]

Therefore
\[
\left(  \sum_{\alpha\in\Gamma_{\Lambda}}|c_{\alpha}(P)|^{\frac{2\dim\Lambda%
	}{1+\dim\Lambda}}\right)  ^{\frac{1+\dim\Lambda}%
	{2\dim\Lambda}}\leq m!\left(  \sum_{\left(  i_{1},\dots,i_{m}\right)
	\in\Lambda}|\hat{P}(e_{i_{1}},\dots,e_{i_{m}})|^{\frac{2\dim
		\Lambda}{1+\dim\Lambda}}\right)  ^{\frac{1+\dim\Lambda}{2\dim
		\Lambda}}%
\]
and
\begin{align*}
\left(  \sum_{\alpha\in\Gamma_{\Lambda}}|c_{\alpha}(P)|^{\frac{2\dim\Lambda%
	}{1+\dim\Lambda}}\right)  ^{\frac{1+\dim\Lambda}%
	{2\dim\Lambda}}  &  \leq\left(  m!\right)  \left(  m\left(  \frac
{2}{\sqrt{\pi}}\right)  ^{m-1}C_{\Lambda}\right)  \left\Vert \hat
{P}\right\Vert \\
&  \leq\left(  m!\right)  \left(  m\left(  \frac{2}{\sqrt{\pi}}\right)
^{m-1}C_{\Lambda}\right)  e^{m}\Vert P\Vert.
\end{align*}
Since we are dealing only with complex scalars, by the Maximum Modulus
Principle, we have
\[
\left(  \sum_{\alpha\in\Gamma_{\Lambda}}|c_{\alpha}(P)|^{2}\right)  ^{\frac{1}{2}}%
\leq\left(  \sum_{|\alpha|=m}|c_{\alpha}(P)|^{2}\right)  ^{\frac{1}{2}}%
\leq\Vert P\Vert.
\]
Since%
\[
\frac{1}{\frac{2m}{m+1}}=\frac{\theta}{\frac{2\dim\Lambda}%
	{1+\dim\Lambda}}+\frac{1-\theta}{2}%
\]
with
\[
\theta=\frac{\dim\Lambda}{m},
\]
by H\"{o}lder's inequality, we have
\begin{align}
\left(  \sum_{\alpha\in\Gamma_{\Lambda}}|c_{\alpha}(P)|^{\frac{2m}{m+1}}\right)
^{\frac{m+1}{2m}}  &  \leq\left[  \left(  \sum_{\alpha\in\Gamma_{\Lambda}}|c_{\alpha
}(P)|^{\frac{2\dim\Lambda}{1+\dim\Lambda}}\right)
^{\frac{1+\dim\Lambda}{2\dim\Lambda}}\right]  ^{\frac
	{\dim\Lambda}{m}}\left[  \left(  \sum_{\alpha\in\Gamma_{\Lambda}}|c_{\alpha
}(P)|^{2}\right)  ^{\frac{1}{2}}\right]  ^{1-\frac{\dim\Lambda}{m}%
}\label{999}\\
&  \leq\left(  \left(  m!\right)  \left(  m\left(  \frac{2}{\sqrt{\pi}%
}\right)  ^{m-1}C_{\Lambda}\right)  e^{m}\Vert P\Vert\right)
^{\frac{\dim\Lambda}{m}}\Vert P\Vert^{1-\frac{\dim\Lambda%
	}{m}}.\nonumber
\end{align}
Thus%
\[
\left(  \sum_{\alpha\in\Gamma_{\Lambda}}|c_{\alpha}(P)|^{\frac{2m}{m+1}}\right)
^{\frac{m+1}{2m}}\leq e^{\dim\Lambda}\left(  C_{\Lambda%
}mm!\right)  ^{\frac{\dim\Lambda}{m}}\left(  \frac{2}{\sqrt{\pi}%
}\right)  ^{\frac{\left(  m-1\right)  \dim\Lambda}{m}}\Vert P\Vert.
\]

\begin{remark}
	If there is a constant $C$ such that $C_{\Lambda}\leq C^{m}$ for all
	$m$, then the constant in the above inequality is dominated asymptotically by%
	\[
	\left(  \frac{2C}{\sqrt{\pi}}\right)  ^{\dim\Lambda}\left(  m\right)
	^{\dim\Lambda}.
	\]
	In fact, by Stirling's Formula we have that%
	\begin{align*}
	e^{\dim\Lambda}\left(  C_{\Lambda}mm!\right)  ^{\frac
		{\dim\Lambda}{m}}\left(  \frac{2}{\sqrt{\pi}}\right)  ^{\frac{\left(
			m-1\right)  \dim\Lambda}{m}} &  \leq\left(  Ce\right)  ^{\dim
		\Lambda}\left(  mm!\right)  ^{\frac{\dim\Lambda}{m}}\left(
	\frac{2}{\sqrt{\pi}}\right)  ^{\frac{\left(  m-1\right)  \dim\Lambda%
		}{m}}\\
	&  \sim\left(  \frac{2Ce}{\sqrt{\pi}}\right)  ^{\dim\Lambda}\left(
	m!\right)  ^{\frac{\dim\Lambda}{m}}\\
	&  \sim\left(  \frac{2Ce}{\sqrt{\pi}}\right)  ^{\dim\Lambda}\left(
	\frac{m^{1+\frac{1}{2m}}}{e}\right)  ^{\dim\Lambda}\\
	&  \sim\left(  \frac{2C}{\sqrt{\pi}}\right)  ^{\dim\Lambda}\left(
	m^{1+\frac{1}{2m}}\right)  ^{\dim\Lambda}\\
	&  \sim\left(  \frac{2C}{\sqrt{\pi}}\right)  ^{\dim\Lambda}\left(
	m\right)  ^{\dim\Lambda}.
	\end{align*}
	
\end{remark}


\begin{thebibliography}{9}                                                                                                %
	
	
	\bibitem {www}N. Albuquerque, G. Ara\'{u}jo, D. N\'{u}\~{n}ez-Alarc\'{o}n, D.
	Pellegrino, and P. Rueda, On summability of multilinear operators and
	applications. Ann. Funct. Anal. \textbf{9} (2018), no. 4, 574--590.
	
	\bibitem {bayart}F. Bayart, Summability of the coefficients of a multilinear
	form, preprint.
	
	\bibitem {bohr}F. Bayart, D. Pellegrino and J. B. Seoane-Sep\'{u}lveda, The
	Bohr radius of the $n$--dimensional polydisk is equivalent to $\sqrt{(\log
		n)/n}$, Adv. Math. \textbf{264} (2014), 726--746.
	
	\bibitem {bh}H. F. Bohnenblust, E. Hille, On the absolute convergence of
	Dirichlet series, Ann. of Math. \textbf{32} (1931), 600--622.
	
	\bibitem {carando}D. Carando, A. Defant, P. Sevilla-Peris, The
	Bohnenblust-Hille inequality combined with an inequality of Helson, Proc.
	Amer. Math. Soc. \textbf{143} (2015), no. 12, 5233--5238.
	
	\bibitem {ann}A. Defant, L. Frerick, J. Ortega-Cerd\'{a}, M. Ouna{\"{\i}}es,
	K. Seip, The Bohnenblust-Hille inequality for homogeneous polynomials is
	hypercontractive, Ann. of Math. (2), \textbf{174} (2011), 485--497.
	
	\bibitem {mais}M. Maia, T. Nogueira, D. Pellegrino, The Bohnenblust-Hille
	inequality for polynomials whose monomials have a uniformly bounded number of
	variables. Integral Equations Operator Theory \textbf{88} (2017), no. 1, 143--149.
	
	
\end{thebibliography}
\end{document}